\documentclass{birkart}
\usepackage{amssymb,amscd}

\DeclareMathOperator{\clos}{clos}
\DeclareMathOperator{\id}{id}
\DeclareMathOperator{\Ran}{Ran}
\DeclareMathOperator{\Ker}{Ker}
\DeclareMathOperator{\Ext}{Ext}
\DeclareMathOperator{\Int}{Int}
\DeclareMathOperator{\Mor}{Mor}
\DeclareMathOperator{\Ob}{Ob}
\DeclareMathOperator{\Sys}{Sys}
\DeclareMathOperator{\Mod}{Mod}
\DeclareMathOperator{\Tfn}{Tfn}
\DeclareMathOperator{\Cfn}{Cfn}

\DeclareMathOperator{\ind}{ind}

\begin{document}

\theoremstyle{plain}
\newtheorem{theorem}{Theorem}[section]
\theoremstyle{definition}
\newtheorem{remark}{Remark}[section]
\theoremstyle{plain}
\newtheorem*{theorema}{Theorem A}
\newtheorem*{theoremb}{Theorem B}
\newtheorem*{theorem_}{Theorem}
\newtheorem{lemma}[theorem]{Lemma}
\newtheorem{prop}[theorem]{Proposition}
\newtheorem{cor}{Corollary}
\renewcommand{\thecor}{}
\newtheorem*{defin}{Definition}
\newcommand{\myproof}{{\noindent\it Proof. }}
\newcommand{\proofend}{$\;\square$\\[3pt]}
\newcommand{\exmpl}{\noindent\textbf{Example.}}
\newcommand{\myrem}{\noindent\textbf{Remark.}}

\title[Conservative curved systems]
{Conservative curved systems and \\ free functional model}

\author[Alexey Tikhonov]{Alexey Tikhonov}

\email{tikhonov@club.cris.net}



\begin{abstract}
We introduce conservative curved systems over multiply connected domains and study relationships of
such systems with related notions of functional model, characteristic function, and transfer
function. In contrast to standard theory for the unit disk, characteristic functions and transfer
functions are essentially different objects. We study possibility to recover the characteristic
function for a given transfer function. As the result we obtain the procedure to construct the
functional model for a given conservative curved system. We employ the functional model to solve
some problems in perturbation theory.
\end{abstract}

\maketitle \thispagestyle{empty}

\setcounter{section}{-1}

\section{Introduction}

It is well known~\cite{NF,Br} that there is a one-to-one correspondence between unitary
colligations or (in terminology of system theory) conservative linear systems
\[
{\mathfrak A}=\left(\begin{array}{cc}T&N\\M&L\end{array}\right) \in{\mathcal L}(H\oplus{\mathfrak
N}, H\oplus{\mathfrak M}),\quad {\mathfrak A}^*{\mathfrak A}^{}=I,\;\,{\mathfrak A}^{}{\mathfrak
A}^*=I
\]
and operator valued functions $\Theta(z)$ of the Schur class
\[
S=\{\Theta\in H^\infty({\mathbb D},{\mathcal L}({\mathfrak N},{\mathfrak M})) :
||\Theta||_\infty\le 1\}\,.
\]
 The mapping defined by the formula $\;\Theta(z)=L^*+zN^*(I-zT^*)^{-1}M^*\,,\,|z|<1$ is
one of the directions of the above mentioned correspondence. The opposite direction is realized via
functional model~\cite{NF,Br}, whose essential ingredients are Hardy spaces $H^2$ and $H^2_-$
(see~\cite{Du}; note that $L^2=H^2 \oplus H^2_-$). These two sides of the theory are equipollent:
both have simple, clear and  independent descriptions and we can easily change a point of view from
unitary colligations to Schur class functions and back.

We shall consider operator valued functions (or rather sets of operator valued functions
$\Theta=(\Theta^{+},\Xi_+,\Xi_-)$) of \textit{weighted Schur classes} $S_{\Xi}$  :
\[
\begin{array}{ll}
S_{\Xi}\;:=&\{\;(\Theta^{+},\Xi_+,\Xi_-)\; :\; \Theta^{+}\in H^\infty(G_+,{\mathcal L}({\mathfrak
N}_+,{\mathfrak N}_-))\,,\\[2pt]
&\quad\forall\;\zeta\in C\;\;\forall\;n\in {\mathfrak N}_+\quad ||\Theta^{+}(\zeta)n||_{-,\zeta}\le
||n||_{+,\zeta}\}\,,
\end{array}\eqno{\rm{(Cfn)}}
\]
where $G_+$ is a multiply-connected domain bounded by a rectifiable Jordan Carleson curve $C$ (we
use also notation $G_-$ for $\Ext C$), $\,\Xi_{\pm}\,$ are weights such that
$\Xi_{\pm},\Xi_{\pm}^{-1}\in L^\infty(C,{\mathcal L}({\mathfrak N}_{\pm}))$, $\Xi_{\pm}(\zeta)\ge
0,\;\zeta\in C\;$, and $\,||n||_{\pm,\zeta}=(\Xi_{\pm}(\zeta)n,n)^{1/2},\;n\in {\mathfrak
N}_{\pm}\,$.

At this moment we ought to look for a suitable generalization of conservative systems (=unitary
colligations), but we find ourself in a position to make a choice what features of unitary
colligations to retain. First, we are going to keep conservatism, that means the characteristic
function of a such system ought to be of weighted Schur class. Here we fork with~\cite{BGK}, where
general transfer functions, rather than conservative ones, is a subject of study. Second, we want
to retain analiticity in both domains $G_+$ and $G_-$. In such a way we reserve to ourself
possibility to exploit techniques that are typical for the boundary values problems (singular
integral operators, the Riemman-Hilbert problem, etc) including connections to the stationary
scattering theory~\cite{Yaf} (e.g., the smooth methods of T.Kato). Thus we will use both Smirnov
spaces $E^2(G_{\pm})$~\cite{Du}, which are analogues of the Hardy spaces $H^2$ and $H^2_-$. The
requirement of analyticity in both domains conflicts with requirement of orthogonality: in general,
the decomposition $L^2(C)=E^2(G_{+})\dot+E^2(G_{-})$ is not orthogonal. The combination
``analyticity in $G_+$ plus orthogonality'' is a mainstream of development in the
multiply-connected case starting with~\cite{AD} and therefore at this point we fork with
traditional way of generalization of Sz.-Nagy-Foia\c{s} theory~\cite{AD,Ba,PF}. Nevertheless, our
requirements are also substantial and descends from applications (see~\cite{T1,Ya,VYa}):
in~\cite{T1} we studied the duality of spectral components for trace class perturbations of a
normal operator with spectrum on a curve; the functional model from~\cite{Ya} goes back
to~\cite{VYa}, which is devoted to spectral analysis of linear neutral functional differential
equations. Moreover, one can regard the present paper as a result of comparison and unification of
the approaches from~\cite{T1,Ya}. If we examine our variant of functional model with more
retrospective point of view, we can observe that it lies on intersection of evolution lines going
back to \cite{NF,Br,AD,MV,Na}.

In order to formulate a main result, we need some background (see for details~\cite{T1,T2,T3}).
Consider pairs $\Pi=(\pi_+,\pi_-)$ of operators $\pi_{\pm}\in {\mathcal L}(L^2(C,{\mathfrak
N_{\pm}}),{\mathcal H})$ such that
\[
\begin{array}{clcl}
(i)_1&(\pi_{\pm}^*\pi_{\pm})z=z(\pi_{\pm}^*\pi_{\pm});&
(i)_2&\pi_{\pm}^*\pi_{\pm}>>0;\\
(ii)_1&(\pi_-^{\dag}\pi_+)z=z(\pi_-^{\dag}\pi_+);&
(ii)_2&P_{-}(\pi_{-}^{\dagger}\pi_{+})P_+=0;\\
(iii)&\Ran\pi_+\vee \Ran\pi_-={\mathcal H}\,,
\end{array}  \eqno{\rm{(Mod)}}
\]
where $P_{\pm}$ are projections such that $\Ran P_{\pm}=E^2(G_{\pm})$ and $\Ker
P_{\pm}=E^2(G_{\mp})\,$. Operators $\;\pi_{\pm}\colon L^2(C,\Xi_{\pm})\to {\mathcal H}\;$ are
isometries if we regard them  as operators acting from weighted $L^2$ with operator valued weights
$\Xi_{\pm}=\pi_{\pm}^*\pi_{\pm}\,$. Note that in this interpretation $\pi_{\pm}^{\dag}$ are the
adjoint operators of $\pi_{\pm}$. We shall say that $\Pi$ is a free \textit{functional model} of
Sz.-Nagy-Foia\c{s} type (see~\cite{NV} for the unit disk case) and write $\Pi\in\Mod$. However, in
contrast to~\cite{NV}, in our axiomatics we use neither a unitary dilation nor orthogonal
compliments.

It can easily be checked that
\[
\Theta=(\pi_{-}^{\dagger}\pi_{+},\pi_{+}^*\pi_{+},\pi_{-}^*\pi_{-})\in S_{\Xi}\,.\eqno{\rm{(MtoC)}}
\]
Therefore we have defined the functor $\mathcal{F}_{cm} :  \Mod\to\Cfn $, where $\Ob(\Cfn)=S_{\Xi}$
and the symbol $\,\Cfn\,$ denotes the category of \textit{characteristic functions}. Conversely,
for a given $\Theta\in S_{\Xi}\,$, it is possible to construct (up to unitarily equivalence) a
functional model $\Pi\in \Mod$ such that
$\Theta=(\pi_{-}^{\dagger}\pi_{+},\pi_{+}^*\pi_{+},\pi_{-}^*\pi_{-})\,$ ,i.e., there exists the
inverse functor $\mathcal{F}_{mc}=\mathcal{F}_{cm}^{-1} : \Cfn\to\Mod$. Note that any
characteristic function $\Theta\in S_{\Xi}\,$ can be transformed to the form
$\Theta_{p}\oplus\Theta_{u}$, where $\Theta_{p}$ is the pure "part" of $\Theta$ and $\Theta_{u}$ is
the unitary "part" of $\Theta$ (see Remark~\ref{rem2}).

We define curved conservative systems in terms of the functional model. Let $\Pi\in\Mod$. Define
the model system
$\,\widehat{\Sigma}=\mathcal{F}_{sm}(\Pi):=(\widehat{T},\widehat{M},\widehat{N},\widehat{\Theta}_u,\widehat{\Xi})\,$,
where\\
\[
\begin{array}{ll}
{\widehat T}\in\mathcal{L}(\mathcal{K}_{\Theta})\,,\;\, & {\widehat T}f:=\mathcal{U}f-\pi_{+}{\widehat M}f\,;\\[2pt]
{\widehat M}\in\mathcal{L}(\mathcal{K}_{\Theta},{\mathfrak N_+})\,,\;\, & {\widehat
M}f:=\displaystyle \frac{1}{2\pi
i}\int_{C} (\pi_{+}^{\dag}f)(z)\,dz \,;\\
{\widehat N}\in\mathcal{L}({\mathfrak N_-},\mathcal{K}_{\Theta})\,,\;\, & {\widehat
N}n:=P_{\Theta}\pi_{-}n \,;\\[6pt]
\widehat{\Theta}_u \;\;\mbox{is the unitary "part" of} &
\widehat{\Theta}=(\pi_-^{\dag}\pi_+,\,\widehat{\Xi})\,;
\\[7pt]
\widehat{\Xi}:=(\pi_+^*\pi_+,\,\pi_-^*\pi_-)\,;
\end{array}
\eqno{\rm{(MtoS)}}
\]
$f\in \mathcal{K}_{\Theta}:=\Ran P_{\Theta}\,$,
$\;P_{\Theta}:=(I-\pi_{+}P_+\pi_{+}^{\dag})(I-\pi_{-}P_-\pi_{-}^{\dag})\,,\;n\in {\mathfrak
N_-}\,$, and the normal operator $\mathcal{U}$ with absolutely continuous spectrum lying on $C$ is
uniquely determined  by conditions $\,\mathcal{U}\pi_{\pm}=\pi_{\pm}z\,$. Note that the model
operators $\;\widehat{T},\widehat{M},\widehat{N}\;$ arise in an natural way and are the simplest in
certain sense (see Remark~\ref{rem3}).

\vskip 2pt A coupling $\Sigma=(T,M,N,\Theta_u,\Xi)$ is called a \textit{curved conservative system}
if there exists a functional model $\Pi$, an operator $W\in
\mathcal{L}(H,\mathcal{K}_{\Theta}\oplus\mathcal{K}_u)$, and a normal operator $\widehat{T}_u\in
\mathcal{L}(\mathcal{K}_u)\,,\,\sigma(\widehat{T}_u)\subset C$ such that $W^{-1}\in
\mathcal{L}(\mathcal{K}_{\Theta}\oplus\mathcal{K}_u,H)$,
\[
(\widehat{T}\oplus\widehat{T}_u)W=WT\,,\quad\widehat{M}W=M\,,\quad\widehat{N}=NW\,,\quad
\widehat{\Theta}_u=\Theta_u\,,\quad\widehat{\Xi}=\Xi\,.\eqno{\rm{(Sys)}}
\]
We will use the symbol $\,\Sys\,$ to denote the category of curved conservative systems. All curved
conservative systems form $\Ob(\Sys)$. We shall mainly consider systems with
$\,\mathcal{K}_u=\{0\}\,$. The relations (MtoS) defines the functor
 $\mathcal{F}_{sm} : \Mod\to\Sys$.

\vskip 2pt For a curved conservative system $\Sigma$ define the \textit{transfer function}
\[
\mathcal{F}_{ts}(\Sigma):=(\Upsilon(z),\,\Theta_u,\,\Xi) \;\mbox{,
where}\quad\Upsilon(z):=M(T-z)^{-1}N\,.\eqno{\rm{(Tfn)+(StoT)}}
\]
The following relation between transfer and characteristic functions is central in our theory
\[
\Upsilon(z)=\left\{%
\begin{array}{ll}
    \Theta^-_+(z)-\Theta^+(z)^{-1}\,, & z\in G_+\cup\rho(T)\,; \\
    -\Theta^-_-(z)\,, & z\in G_-\,,
\end{array}%
\right.\eqno{\rm{(CtoT)}}
\]
where the operator valued functions $\Theta^-_{\pm}(z)$ are defined by
\[
\begin{array}{l}
    \Theta^-_{\pm}(z):=(P_{\pm}\Theta^-n)(z),\;\;z\in G_{\pm},\;\;n\in \mathfrak{N}_-\,; \\[3pt]
    \Theta^-(\zeta):=(\pi_+^{\dag}\pi_-)(\zeta)=\Xi_+(\zeta)^{-1}\Theta^+(\zeta)^*\Xi_-(\zeta),\;\;\zeta\in C\,.
\end{array}
\]
The relation (CtoT) gives the expression for transfer functions in terms of a characteristic
functions and uses both metric and analytic properties of a characteristic function. The relations
(StoT) and (CtoT) define the functors $\mathcal{F}_{ts} : \Sys\to\Tfn$ and $\mathcal{F}_{ts} :
\Cfn\to\Tfn\;$, respectively. Thus we have arrived to the following commutative diagram
\[
\begin{CD}
  \Mod @>\mathcal{F}_{cm}>> \Cfn\\
  @V\mathcal{F}_{sm}VV  @VV\mathcal{F}_{tc}V\\
  \Sys @>>\mathcal{F}_{ts}> \Tfn
\end{CD}\eqno{\rm{(C\&F)}}
\]

We have defined the notion of conservative curved system.  Linear similarity (instead of unitarily
equivalence for unitary colligations) is a natural kind of equivalence for conservative curved
systems and duality is a substitute for orthogonality. Another distinctive feature of conservative
curved systems is that we distinguish notions of characteristic and transfer function. In Section 1
we give a more detailed survey of categories $\,\Mod,\Cfn,\Sys,\Tfn\,$ and related functors. The
fundamental relation $\,\rm{(\Phi\&F)}\,$ from that Section explains a functorial character of the
transformation $\,\mathcal{F}_{cm},\,\mathcal{F}_{sm},\,\mathcal{F}_{tc},\,\mathcal{F}_{ts}\,$.

As we can now see, characteristic functions and conservative curved systems are not on equal terms:
first of them plays leading role because the definition of conservative curved system depends on
the functional model, which, in turn, is uniquely determined by the characteristic function. But,
surprisingly, the conservative curved systems is a comparatively autonomous notion and one of the
aims of this paper is to ``measure'' a degree of this autonomy.

Examining the diagram (C\&F) we see that it is an interesting problem to recover a functional model
$\Pi\in\Mod$ for a given conservative curved systems $\Sigma\in\Sys$. We have
$\,\mathcal{F}_{sm}^{-1}=\mathcal{F}_{cm}^{-1}\circ\mathcal{F}_{tc}^{-1}\circ\mathcal{F}_{ts}$.
Since the functor $\mathcal{F}_{cm}$ is evidently invertible, the main problem is to invert
$\mathcal{F}_{tc}$. In fact, the projections $P_{\pm}$ are singular integral operators and we can
regard (CtoT) as a system of singular integral equations with unknown $\Theta^+$ and given
$\Upsilon,\Xi_{\pm}$\,. The main result of Section~2 (and of the paper on the whole) is the
uniqueness of a solution of this system.
\begin{theorema}
If $\;\Theta_1,\,\Theta_2\in\Cfn\,$ and
$\,\mathcal{F}_{ct}(\Theta_1)=\mathcal{F}_{ct}(\Theta_2)\in\mathcal{N}(G_+\cup
G_-,\mathcal{L}(\mathfrak{N}_{-},\mathfrak{N}_{+}))\,$, then $\;\Theta_1=\Theta_2\,$.
\end{theorema}
\noindent In other words, if the transfer function $\Upsilon$ is a function of the Nevanlinna class
and the weights $\Xi_{\pm}$ are fixed, then the pure "part" $\Theta_p^+$ of the characteristic
function is uniquely determined. Here $\,\mathcal{N}(G_+\cup
G_-,\mathcal{L}(\mathfrak{N}_{-},\mathfrak{N}_{+}))$ is the Nevanlinna class of operator valued
functions that admit representation of the form $\Upsilon(z)=1/\delta(z)\,\Omega(z)\,$, where
$\delta\in H^{\infty}(G_+\cup G_-)$ and $\Omega\in H^{\infty}(G_+\cup
G_-,\mathcal{L}(\mathfrak{N}_{-},\mathfrak{N}_{+}))$.

We also describe a procedure how to find a characteristic function for given transfer function.
Knowing the characteristic function we can construct the functional model $\Pi$ and obtain explicit
formulas for conversion of the simple system $\Sigma$ into the model system $\widehat{\Sigma}$. As
a result we give purely intrinsic description of simple conservative curved systems with the
transfer function of the Nevanlinna class. Note that conservative curved systems arise for certain
problems of operator theory and we reap the benefit of functional model when we are able to
determine that some set of operators $(T,M.N)$ is a curved conservative system~\cite{T1,Ya}.

In Section 3 we employ the functional model to establish duality of spectral components of trace
class perturbations of normal operators with spectrum on a curve. The statement is the same as
in~\cite{T1}, but now we extend that assertion to the multiply connected case. One of our purposes
is to demonstrate how and why the functional model for conservative curved systems can be
successfully used in perturbation theory.

\vskip 2pt The author is grateful to D.Yakubovich for interesting and useful information and J.Ball
for stimulating questions.

\section{Categories and functors}

In this Section we survey (mainly without proofs) basic facts concerning the categories
$\,\Mod,\Cfn,\Sys,\Tfn\,$ and related functors (see~\cite{T1,T2,T3} for details). We shall use the
language of categories and functors for short and for more systematic character of it (anyway we
need some notation for our classes of objects and their transformations). We deal with the
categories of \textit{functional models} $\,\Mod\,$, \textit{characteristic functions} $\,\Cfn\,$,
\textit{conservative curved systems} $\,\Sys\,$, and \textit{transfer functions} $\,\Tfn\,$. Let
$\,X\in\{\Mod,\Cfn,\Sys,\Tfn\}\,$, e.g., $X=\Mod$. Define objects of category $X$
\[
\Ob(X):=\{ O : O \;\,\mbox{satisfies the condition}\;(X)\,\}\,.\eqno{(\Ob)}
\]
We regard two models $\,\Pi_1,\Pi_2\in\Ob(\Mod)\,$ as equal if there exists an unitary operator
$\,V : \mathcal{H}_1\to\mathcal{H}_2\,$ such that $\,\pi_{2\pm}=V\pi_{1\pm}\,$.\\[2pt]
We regard two systems $\,\Sigma_1,\Sigma_2\in\Ob(\Sys)\,$ as equal if there exists an invertible
operator $W\in\mathcal{L}(H_1,H_2)$ such that
\[
T_2W=WT_1\,,\quad M_2W=M_1\,,\quad N_2=N_1W\,,\quad\Theta_{2u}=\Theta_{1u}
\,,\quad\Xi_{2\pm}=\Xi_{1\pm}\;.
\]
\begin{remark}\label{rem1} In this paper we deal with generalization of functional models
from~\cite{T1,Ya}. In~\cite{T1} we made use of models with
$\,\Xi_{\pm}=\pi_{\pm}^*\pi_{\pm}=\delta_{\pm}(\zeta)I\,,\;\zeta\in C\,$, where
$\delta_{\pm}(\zeta)$ were scalar bounded positive functions. In~\cite{Ya} D.Yakubovich developed
his theory using models such that $\,\Xi_+=I\,$ and $\,\Xi_-=(\Theta^{+}\Theta^{+*})^{-1}\,$,
though he did not introduce the weights $\,\Xi_{\pm}\,$ in explicit form.
\end{remark}
\begin{remark}\label{rem2}
A characteristic function $\,\Theta\in\Ob(\Cfn)\,$ (for short we shall often use
$\,\Theta\in\Cfn\,$) is called $\Xi$-pure if  the (multiple valued character-automorphic) operator
valued function $\,\Theta_0^+=\chi_-\Theta^+\chi_+^{-1}\,$ satisfies the condition $\,\forall z\in
G_+\;\forall n\in\mathfrak{N}_+,\;n\ne0\,$\\
$\,||\Theta_0^+(z)n||_{\mathfrak{N}_-}<||n||_{\mathfrak{N}_+}\,$, where $\,\chi_{\pm}\,$ are outer
(character-automorphic) operator valued functions such that $\,\chi_{\pm}^*\chi_{\pm}=\Xi_{\pm}\,$.
A characteristic function $\,\Theta\in\Cfn\,$ is called a $\Xi$-unitary constant if
$\,\Theta^-=(\Theta^+)^{-1}\in H^{\infty}(G_+,\mathcal{L}(\mathfrak{N}_{-},\mathfrak{N}_{+}))$. For
any $\,\Theta\in\Cfn\,$, there exist $\eta_{\pm}$ such that $\,\eta_{\pm},\,\eta_{\pm}^{-1}\in
H^{\infty}(G_+,\mathcal{L}(\mathfrak{N}_{\pm}))\,$ and the characteristic function
$\,\Theta_{pu}=(\eta_-^{-1}\Theta^+\eta_+,\,\eta_+^*\Xi_+\eta_+,\,\eta_-^*\Xi_-\eta_-)\in\Cfn\,$
can be represented in the form $\,\Theta_{pu}=\Theta_{p}\oplus\Theta_{u}\,$, where $\,\Theta_{p}\,$
is  $\Xi$-pure and $\,\Theta_{u}\,$ is a $\Xi$-unitary constant. If
$\,\Theta_{pu}'=\Theta_{p}'\oplus\Theta_{u}'\,$ is another such representation for some
$\eta_{\pm}'$, then the operator valued functions $\psi_{\pm}=\eta_{\pm}^{-1}\eta_{\pm}'$ admit the
decompositions $\psi_{\pm}=\psi_{\pm}^p\oplus\psi_{\pm}^u :
{\mathfrak{N}_{\pm}^p}'\oplus{\mathfrak{N}_{\pm}^u}'\to\mathfrak{N}_{\pm}^p\oplus\mathfrak{N}_{\pm}^u$.
This enable us to consider characteristic functions with equal $\Xi$-pure and/or $\Xi$-unitary
"parts".
\end{remark}
\begin{remark}\label{rem3}
For a functional model $\Pi\in\Mod$, one can define two subspaces
$\,\mathcal{D}_+:=\Ran\,\pi_+P_+\pi_+^{\dag}\,$ and
$\,\mathcal{H}_+:=\Ran\,(I-\pi_-P_-\pi_-^{\dag})\,$. It is easily shown that
$\,\mathcal{D}_+\subset\mathcal{H}_+\,$ and $\,\forall\,z\in
G_+\quad\mathcal{D}_+\subset(\mathcal{U}-z)^{-1}\mathcal{D}_+\,,
\quad\mathcal{H}_+\subset(\mathcal{U}-z)^{-1}\mathcal{H}_+\,$. In particular,
$\mathcal{D}_+\subset\mathcal{U}\mathcal{D}_+\,,\;\mathcal{H}_+\subset\mathcal{U}\mathcal{H}_+\,$.
Whence we see that the operator $\,\mathcal{U}\,$ is a normal dilation of the operator
$\,\widehat{T}=P_{\Theta}\mathcal{U}|\mathcal{K}_{\Theta}\,$. After straightforward computation we
get $\,\widehat{T}f=\mathcal{U}f-\pi_+\widehat{M}f\,$. Taking into account that $\Ran
P_{\pm}=E^2(G_{\pm})$ and $\Ker P_{\pm}=E^2(G_{\mp})$, we obtain
$\,\widehat{M}f=(\pi_+^{\dag}\mathcal{U}f)(\infty),\;\;\widehat{M}
 : \mathcal{K}_{\Theta}\to\mathfrak{N}_+\,$ and
\[
(\widehat{T}-z)^{-1}f=(\mathcal{U}-z)^{-1}(f-\pi_+n)\,, \quad n=\left\{\begin{array}{ll}
\Theta^+(z)^{-1}(\pi_-^{\dag}f)(z)&,\,z\in G_+\cap\rho(\widehat{T})\\
(\pi_+^{\dag}f)(z)&,\,z\in G_-\;
\end{array}.\right.
\]
Assuming the similar conditions for the dual model (see below) and using the duality $\,({\widehat
N}n,g)=(n,\widehat{M}_*g)\,,\;g\in\mathcal{K}_{*\Theta},\,n\in\mathfrak{N}_-\,$, we obtain
$\,{\widehat N}n=P_{\Theta}\pi_{-}n\,$ . If we hold the condition $\Ran P_{+}=E^2(G_{+})$ and drop
$\Ker P_{+}=E^2(G_{-})$, we have merely $\,\widehat{M} : \mathcal{K}_{\Theta}\to
L^2(C,\mathfrak{N}_+)\,$ and nothing more. Note that it is possible to consider the functional
model under those more weak assumptions, but then we lose one of advantages of our functional
model, namely, its ability for automatic, easy and effective computation (see Section~3).
\end{remark}
\begin{remark}\label{rem4}
It can be shown that
\[
\hat{r}_{nz}= (\widehat{T}-z)^{-1}\widehat{N} n=P_{\Theta}\pi_{\pm}\frac{m}{\zeta-z}\,,\quad
m=\left\{\begin{array}{rl}
-\Theta^+(z)^{-1}n&,\,z\in G_+\cap\rho(\widehat{T})\\
n&,\,z\in G_-\;
\end{array}.\right.
\]
If $\,\rho(\widehat{T})\cap G_+\ne\emptyset\;$, then
$\,\vee_{z\in\rho(\widehat{T})}\hat{r}_{nz}=\mathcal{K}_{\Theta}\,$. In this connection, we call a
system $\,\Sigma\in\Sys\,$ simple if $\,\rho({T})\cap G_+\ne\emptyset\;$ and
$\,\vee_{z\in\rho({T})}{r}_{nz}=H\,$, where
\[
r_{nz}:= (T-z)^{-1}N n\,,\quad n\in \mathfrak{N}_-\,,\quad z\in\rho(T)\,.\qquad \eqno{\rm{(RK)}}
\]
\end{remark}
\vskip 6pt We postpone to introduce morphisms for our categories and at first consider functors
$\,\mathcal{F}_{YX} : \Ob(X)\to\Ob(Y)\,$, where $\,X,Y\in\{\Mod,\Cfn,\Sys,\Tfn\}\,$. In the
Introduction we have defined the functors
$\,\mathcal{F}_{cm},\,\mathcal{F}_{sm},\,\mathcal{F}_{tc},\,\mathcal{F}_{ts}\,$. By means of these
functors we can translate relationships from language of one category into language of a parallel
category. For instance, it is possible to give a complete description of the spectrum of a model
operator $\,\widehat{T}\,$ in terms of its characteristic function. Another example is the
existence of the one-to-one correspondence between regular factorizations of a characteristic
function and invariant subspaces of the operator~$\,\widehat{T}\,$.

Together with functors $\,\mathcal{F}_{YX}\,$ we consider transformations $\,\Phi_{\eta}^{X} :
\Ob(X)\to\Ob(X)\,$, where $\,\eta=(\varphi,\eta_+,\eta_-)\in CM_{\eta}\,$. The latter means
$\,\varphi : G_{1+}\to G_{2+}\,$ is a conformal mapping and $\,\eta_{\pm},\eta_{\pm}^{-1}\in
H^{\infty}(G_{2+},\mathcal{L}(\mathfrak{N}_{\pm}))\,$.

For the category $\,\Mod\,$, we set $\,
\Pi_2=\Phi_{\eta}^{\Mod}(\Pi_1):=(\pi_{1+}C_{\varphi}\eta_+,\,\pi_{1-}C_{\varphi}\eta_-)\,$, where
the operator $\,(C_{\varphi}f(\cdot))(\zeta):=\sqrt{\varphi'(\zeta)}f(\varphi(\zeta))\,,\;\zeta\in
C_2\,,\;f\in L^2(C_2,\mathfrak{N})\,$ is unitary.

For the category $\,\Cfn\,$ we set $\,
\Theta_2=\Phi_{\eta}^{\Cfn}(\Theta_1):=(\,\eta_-^{-1}(\Theta_1^+\circ\varphi^{-1})\,\eta_+,
\,\eta_+^*(\Xi_{1+}\circ\varphi^{-1})\,\eta_+,\,\eta_-^*(\Xi_{1-}\circ\varphi^{-1})\,\eta_-)\,$. In
fact, we have made use of this transformation in Remark~\ref{rem2}.

For the category $\,\Sys\,$ we set $\,
\Sigma_2=\Phi_{\eta}^{\Sys}(\Sigma_1):=(\,T_2,\,M_2,\,N_2,\,\Theta_{2u},\,\Xi_2)\,$, where
$\,T_2=\varphi(T_1)\,$,
\[
\begin{array}{c}
\begin{array}{l}
M_2f=\displaystyle-\frac{1}{2\pi
i}\int\limits_{C_1}\sqrt{\varphi'(\zeta)}\;\eta_+(\varphi(\zeta))^{-1}[M_1(T_1-\cdot)^{-1}f]_-(\zeta)\,d\zeta\,,\\[8pt]
N_{2}^*g=\displaystyle-\frac{1}{2\pi
i}\int\limits_{C_{1*}}\sqrt{\bar{\varphi}'(\bar{\zeta)}}\;\eta_-({\varphi}(\bar{\zeta}))^{*}\,[N_1^*(T_1^*-\cdot)^{-1}g]_-(\zeta)\,d\zeta\,,
\end{array}%
\end{array}
\]
and $\,\Xi_{2\pm}=\eta_{\pm}^*(\Xi_{1\pm}\circ\varphi^{-1})\,\eta_{\pm}\,$.

\vskip 5pt The following important identity indicates the connection between functors and the
transformations $\,\Phi_{\eta}^{X}\,$.
\[
\mathcal{F}_{YX}\circ\Phi_{\eta}^{X}=\Phi_{\eta}^{Y}\circ\mathcal{F}_{YX}\,.\qquad\eqno{\rm{(\Phi\&F)}}
\]
This relation  can be verifyed for $\,X=\Mod,\,Y=\Cfn\,$ and $\,X=\Mod ,\,Y=\Sys \,$. In
particular, $\,\mathcal{F}_{sm}\Phi_{\eta}^{\Mod}=\Phi_{\eta}^{\Sys}\mathcal{F}_{sm}\,$ results
from the identities
\[
\widehat{T}_2P_{2\Theta}=P_{2\Theta}T_2,\quad \widehat{M}_2P_{2\Theta}=M_2,\quad
\widehat{N}_2=P_{2\Theta}N_2,\quad\\[3pt]
\]
where $\,\Pi_2=\Phi_{\eta}^{\Mod}(\Pi_1)\,$,
$\;(T_2,\,M_2,\,N_2)=\Phi_{\eta}^{\Sys}(\widehat{T}_1,\,\widehat{M}_1,\,\widehat{N}_1)\,$,
$\,(\widehat{T}_1,\,\widehat{M}_1,\,\widehat{N}_1)=\mathcal{F}_{sm}(\Pi_1)\,$, and
$\,(\widehat{T}_2,\,\widehat{M}_2,\,\widehat{N}_2)=\mathcal{F}_{sm}(\Pi_2)\,$.

\begin{remark}\label{rem5} The linear similarity
$\,\widehat{T}_2P_{2\Theta}=P_{2\Theta}\varphi(\widehat{T}_1)\,$ explains why we deal with
equivalence classes of conservative curved systems rather than merely with systems. Note that only
at first sight the linear similarity arises because of weights $\,\Xi_{\pm}\,$. Actually, the main
reason is the following observation. Let $\,H=H_1\oplus H_2'=H_1\oplus H_2''\,$ and the subspace
$\,H_1\,$ be an invariant under an operator~$\,T\,$. Then in the corresponding triangular
representations
\[
T=\left(%
\begin{array}{cc}
  T_1 & * \\
  0 & T_2' \\
\end{array}%
\right)\quad \mbox{and}\quad T=\left(%
\begin{array}{cc}
  T_1 & * \\
  0 & T_2'' \\
\end{array}%
\right)
\]
the operators $\,T_2'\,$ and $\,T_2''\,$ are similar.
\end{remark}
\vskip 1pt  Conversely, we define the transformation $\,\Phi_{\eta}^{\Tfn}\,$ using
$\,\rm{(\Phi\&F)}\,$ with $\,X=\Cfn\,$ and $\,Y=\Tfn\,$ as a definition. This transformation is
well defined because of the implication $\,\mathcal{F}_{tc}(\Theta_1)=\mathcal{F}_{tc}(\Theta_2)\,$
$\Longrightarrow$
$\,\mathcal{F}_{tc}(\Phi_{\eta}^{\Cfn}(\Theta_1))=\mathcal{F}_{tc}(\Phi_{\eta}^{\Cfn}(\Theta_2))\,$,
which holds for $\,\varphi'\in\cup_{p>1}L^p(C_1)\,$.\,\, There are the explicit formulas
\[
\begin{array}{lll}
\Upsilon_2(z)=\upsilon_2^{\sim}(z)\,,\quad &
\;\upsilon_2(\bar{z})\,n=F_{\eta_{-}^{\sim},z}((\upsilon_1^{\sim}(\cdot)\,n)_-)(\bar{z}) \,,\\[3pt]
& \;\upsilon_1(z)\,n=F_{\eta_{+}^{-1},\varphi}((\Upsilon_1(\cdot)\,n)_-)(z)\,,\quad &z\in G_{2-}\,,
\end{array}
\]
where $\,(\Upsilon_1(\cdot)\,n)_{\pm}(\zeta)\,$ are the boundary values of the vector valued
functions $\,\Upsilon_1(z)\,n\,$ from the domains $\,G_{\pm}\,$,
$\,F_{\eta,\varphi}(u)(z):=[P_-(\eta(\zeta)(u\circ\varphi^{-1})(\zeta))](z)\,$, and
$\,A^{\sim}(z):=A(\bar{z})^*\,$. Taking into account that $P_-$ is actually a singular integral
operator, we see that the transformation $\,F_{\eta,\varphi}\,$ is defined for
$\,\varphi'\in\cup_{p>1}L^p(C_1)\,$ and $\,u\in\cap_{q>1}L^q(C_1,\mathfrak{N})\,$.
\begin{remark}\label{rem6}
The transformation $\,F_{\eta,\varphi}\,$ is an analogue of the well known in complex analysis the
Faber transformation~\cite{Su,Ga}. The classic Faber transformation for a simple connected domain
$G_+$ is the transformation $\,F_{\eta,\varphi}\,$ such that $\,\eta\equiv I\,$, $\,\varphi\,$ is a
conformal map of the external domains normed at the infinity, and the projection $P_-$ is replaced
by $P_+$.
\end{remark}
If $z\in G_{2+}$, we need some extra assumptions to write an explicit formula. It suffices to admit
the existence of boundary values of operator valued function $\,\Theta_1^+(z)^{-1}\,$. In this case
$\,\Upsilon_2(z)\,$ is the analytic continuation of
\[
\Upsilon_2(\varphi(\zeta))_+=\Upsilon_2(\varphi(\zeta))_- +
\eta_+^{-1}(\varphi(\zeta))\,(\Upsilon_1(\zeta)_+ -
\Upsilon_1(\zeta)_-)\,\eta_-(\varphi(\zeta))\,,\quad \zeta\in C_1\,
\]
into the domain $\,G_{2+}\,$.

\vskip 2pt  The transformations $\,\Phi_{\eta}^{X}\,$ is useful if we want to transfer an object
$\,O\in\Ob(X)\,$ to a more simple form, e.g., we can pass to a circular domain or get rid of
weights on some connected component of the boundary. For simple connected domains, in such a way we
can reduce our problems to the case of the unit circle and get rid of the weights $\,\Xi_{\pm}\,$
completely.
\begin{remark}\label{rem7}
In terms of the transformations $\,\Phi_{\eta}^{\Tfn}\,$ one can define $\Xi$-pure transfer
functions. Under the assumption $\,\rho(T)\cap G_+\ne\emptyset\,$ we shall call a transfer function
$\,\Upsilon\in\Ob(\Tfn)\,$ $\Xi$-pure if $\, \forall\,\eta\in CM_{\eta}\quad
\Phi_{\eta}^{\Tfn}(\Upsilon)(z)\,n\equiv0\; \Rightarrow \; n=0\,$. This definition agrees with the
corresponding definition for characteristic functions. Therefore there exists a decomposition
$\,\Upsilon=\Upsilon_p\oplus\Upsilon_u\,$ in the same sense as in Remark~\ref{rem2} (moreover,
$\Upsilon_u\equiv 0$). Note that even in the simplest case of unitary colligations an unitary
matrix $\,\left(\begin{array}{cc}T&N\\M&L\end{array}\right)\,$ is recovered from the triple
$\,(T,M,N)\,$ only up to the unitary part $L_u$ of the operator $L$, where $\,L_u:=L|\Ker(I-L^*L) :
\Ker(I-L^*L) \stackrel{onto}{\longrightarrow} \Ker(I-LL^*)\,$.
\end{remark}
The following properties are basic for the transformations $\,\Phi_{\eta}^{X}\,$:
\[
\Phi_{\id}^{X}=\id_X\,,\qquad
\Phi_{\eta_{32}}^{X}\circ\Phi_{\eta_{21}}^{X}=\Phi_{\eta_{32}\circ\eta_{21}}^{X}\,,\qquad
\eqno{\rm{(\Phi_{\eta})}}
\]
where
$\,\eta_{32}\circ\eta_{21}:=(\varphi_{32}\circ\varphi_{21},\,\eta_{3+}\cdot(\eta_{2+}\circ\varphi_{32}^{-1})
,\,\eta_{3-}\cdot(\eta_{2-}\circ\varphi_{32}^{-1}))\,$. For the category $\,\Mod\,$, these
properties are obvious. For the rest categories they follow from $\,\rm{(\Phi\&F)}\,$.
\begin{remark}\label{rem8}
In the case of the categories $\Cfn$ and $\Sys$ it is possible to give direct proofs, which do not
employ $\,\rm{(\Phi\&F)}\,$.
\end{remark}
At this point we have already prepared to introduce morphisms for our categories. Let
$\,O_1,O_2\in\Ob(X)\,$. We shall say that a triple $m_{O_1 O_2}^X=(O_1,O_2,\eta)$ is a morphism in
the category $\,X\,$ if there exists $\,\eta\in CM_{\eta}\,$ such that
$\,O_2=\Phi_{\eta}^{X}(O_1)\,$. In other words, $\, \Mor(O_1,O_2):=\{ m_{O_1 O_2}^X :
\exists\,\eta\in CM_{\eta}\quad O_2=\Phi_{\eta}^{X}(O_1) \}\,$. The identity $\,m_{O_1 O_3}=m_{O_2
O_3}\circ m_{O_1 O_2}\,$ follows from $\,\rm{(\Phi_{\eta})}\,$. Taking into account
$\,\rm{(\Phi\&F)}\,$, it is easily shown that the transformation $\,\mathcal{F}_{YX} :
\Ob(X)\to\Ob(Y)\,$ can be extend to morphisms and thus $\,\mathcal{F}_{YX}\,$ is a covariant
functor.

Now we pass to the important and useful notion of duality. In language of categories and functors
this means that there exists the transformation $\,{\mathcal{F}_{*X} : \Ob(X)\to\Ob(X)}\,$ such
that $\;\mathcal{F}_{*X}\circ\Phi_{\eta}^{X}=\Phi_{\eta_{-*}}^{X}\circ\mathcal{F}_{*X}\,$ and
$\,\mathcal{F}_{*X}^2=\id_X\,$, where
$\,\eta_{-*}=(\varphi^{\sim},\eta_-^{\sim-1},\eta_+^{\sim-1})\,$. This transformation can be
extended to morphisms and therefore $\,\mathcal{F}_{*X}\,$ is a contravariant functor in the
category $X$. In particular, we have
\[
\begin{array}{l}
\Theta_*=\mathcal{F}_{*c}(\Theta):=(\Theta^{+\sim},\Xi_*)\,,\quad
\Xi_{*\pm}=\Xi_{\mp}^{\sim-1}\,;\\[5pt]
\Pi_*=\mathcal{F}_{*m}(\Pi)\; : \;\;
(f,\pi_{*\mp}v)_{\mathcal{H}}=<\pi_{\pm}^{\dag}f,v>_C\,,
\;f\in\mathcal{H}\,,\;v\in L^2(\bar{C},\mathfrak{N}_{\mp})\,;\\[5pt]
\Sigma_*=\mathcal{F}_{*s}(\Sigma):=(T^*, N^*, M^*, \Theta_u^{\sim}, \Xi_*)\,;\\[5pt]
\Upsilon_*=\mathcal{F}_{*t}(\Upsilon):=(\Upsilon^{\sim}, \Theta_u^{\sim}, \Xi_*)\,,
\end{array}
\]
where $\,<u,v>_C:=\frac{1}{2\pi i}\int_C\;(u(z),v({\bar z)})_{\mathfrak N}\,dz\;,\;u\in
L^2(C,{\mathfrak N}),\;v\in L^2({\bar C},{\mathfrak N})\,$. For the functors
$\,\mathcal{F}_{cm},\,\mathcal{F}_{sm},\,\mathcal{F}_{tc}\,$, the identity
$\,\mathcal{F}_{*Y}\mathcal{F}_{YX}=\mathcal{F}_{YX}\mathcal{F}_{*X}\,$ can easily be checked.

\section{From system to model}

Our aim is to construct the functional model for a given conservative curved system. The main step
on this way is Theorem A, which, in fact, shows that under certain conditions all our categories
$\,\Mod,\Cfn,\Sys,\Tfn\,$ are diverse views of the only one entity. We divide the proof into
several parts, which we arrange as separate assertions.

The main hypothesis of Theorem A is  $\,\Upsilon\in\mathcal{N}(G_+\cup
G_-,\mathcal{L}(\mathfrak{N}_{-},\mathfrak{N}_{+}))\,$. Note that for this it suffices to assume
$\;T-U\in\mathfrak{S}_1,\,U^*U=UU^*,\,\sigma(U)\subset C,\,\sigma_c(T)\subset
C\,,\,M,N\in\mathfrak{S}_2\,$, and $\,C^{1+\varepsilon}$ smooth of the curve $C$
(see~\cite{T1,T4}). But instead of $\,\Upsilon\in\mathcal{N}(G_+\cup
G_-,\mathcal{L}(\mathfrak{N}_{-},\mathfrak{N}_{+}))\,$ we shall consider the more weaker condition,
which consists in the existence of boundary values for an operator valued
function~$\,\Theta^+(z)^{-1}\,$.

\begin{prop}\label{delta}
Suppose $\Theta\in\Cfn$  and $\Theta^+(\zeta)^{-1}$ possesses boundary values a.e. on $C$.
 Then the operator valued function
$\Delta^+(\zeta):=(I-\Theta^+(\zeta)\Theta^-(\zeta))^{1/2}$ can be expressed in terms of the
transfer function $\Upsilon(z)$.
\end{prop}
\myproof The identity
$\,\Theta^-(\zeta)-\Theta^+(\zeta)^{-1}=\Upsilon(\zeta)_+-\Upsilon(\zeta)_-\,,\;\zeta\in C\,$
follows from (CtoT). Then, taking into account that $\,\Theta^-(\zeta)\,$ is adjoint to
$\,\Theta^+(\zeta) : \mathfrak{N}_{+,\zeta} \to \mathfrak{N}_{-,\zeta}\,$ (see (Cfn)), it suffices
to prove the following Lemma.
\begin{lemma}\label{l_delta}
Suppose $\,L\in\mathcal{L}(H_1,H_2)\,$, $\,L^{-1}\in\mathcal{L}(H_2,H_1)\,$, and $\,||L||\le 1\,$.
Let $\,L=|L|U\,$ be the polar decomposition of $L$. Then $\,|L|=\psi(B^*B)\,$, where
$\,B=L^{-1}-L^*\,$ and $\,U^*|\Ran (I-|L|^2)=B(|L|^{-1}-|L|)^{-1}|\Ran (I-|L|^2)\,$.
\end{lemma}
\myproof  $\,\psi(z)=\frac{1}{2}\,\sqrt{2+z-\sqrt{z^2+4z}}\,$\,.\, $\;\square$ \proofend

Let $\,\Pi\in\Mod\,$ and $\,\tau_{+}:=((\Delta^{+})^{-1}
(\pi_{-}^{\dag}-\Theta^{+}\pi_{+}^{\dag}))^{\dag}$. Define the unitary operator $\,W_{NF} :
\mathcal{H} \to \mathcal{H}_{NF}\,$ by the rule
$\,W_{NF}f=(f_{\pi},f_{\tau}):=(\pi_+^{\dag}f,\tau_+^{\dag}f)\,$, where
$\,\mathcal{H}_{NF}:=L^2(C,\Xi_+)\oplus\clos L^2(C,\Xi_-)\,$ is endowed with the inner product
$\,(W_{NF}f,W_{NF}g):=(f_{\pi},g_{\pi})_{L^2(C,\Xi_+)}+(f_{\tau},g_{\tau})_{L^2(C,\Xi_-)}\,$.
Obviously, $\,W_{NF}^{-1}(f_{\pi},f_{\tau})=\pi_+f_{\pi}+\tau_+f_{\tau}\,$. Define also
$\;\mathcal{K}_{\Theta NF}:=W_{NF}\mathcal{K}_{\Theta}\,$.
\begin{prop}\label{rk}
Suppose $\Theta_1,\Theta_2\in\Cfn$;\, $\Theta_1^+(\zeta)^{-1},\Theta_2^+(\zeta)^{-1}$ possesses
boundary values a.e. on $C$, and
$\,\Upsilon=\mathcal{F}_{ct}(\Theta_1)=\mathcal{F}_{ct}(\Theta_2)\,$. Then $\,\mathcal{K}_{1\Theta
NF}=\mathcal{K}_{2\Theta NF}\,$.
\end{prop}
\myproof Take the vectors $\,\hat{r}_{nz}\,$ as in Remark~\ref{rem4}. By direct calculations, we
have
\[
\pi_+^{\dag}\hat{r}_{nz}=\frac{(\Upsilon n)(\zeta)_--\Upsilon(z)\,n}{\zeta-z}\,,\qquad
\tau_+^{\dag}\hat{r}_{nz}=\frac{-\Delta^+(\zeta)\,n}{\zeta-z}\,.
\]
Since $\,\vee_{z\in\rho(\widehat{T})}\hat{r}_{nz}=\mathcal{K}_{\Theta NF}\,$, it remains to make
use of Prop.~\ref{delta}. \proofend
\begin{prop}\label{mod_space}
Let $\,\Theta_1=(\Theta_1^+,\Xi)\,$ and $\,\Theta_2=(\Theta_2^+,\Xi)\,$. Let $\,\chi_{\pm}\,$ be
outer (character-automorphic) operator valued functions such that
$\,\chi_{\pm}^*\chi_{\pm}=\Xi_{\pm}\,$. Suppose $\,\mathcal{K}_{1\Theta NF}=\mathcal{K}_{2\Theta
NF}\,$. Then there exists an unitary operator $\,U\in\mathcal{L}(\mathfrak{N}_-)\,$ such that
$\,\chi_-\Theta_2^+\chi_+^{-1}=U\chi_-\Theta_1^+\chi_+^{-1}\,$.
\end{prop}
\myproof It is obvious that $\,W_{kNF}\,\pi_{k+}=I,\,k=1,2\,$. Hence,
$\,W_{1NF}\pi_{1+}=W_{2NF}\pi_{2+}\,$. On the other hand, the subspace
$\,\mathcal{H}_{NF+}:=\mathcal{K}_{1\Theta NF}\dot+E^2(G_+,\mathfrak{N}_+)\,$ are invariant under
the operator of multiplication by the independent variable $z$. Therefore $\,z|\mathcal{H}_{NF+}\,$
is a subnormal operator with the normal spectrum on the curve $C$. By~\cite{AD}, there exists the
Wold-Kolmogorov decomposition of the space $\,\mathcal{H}_{NF+}\,$ with respect to the operator
$z$. On account of the uniqueness of this decomposition, we get for the pure part
$\,W_{1NF}\pi_{1-}E^2(G_+,\mathfrak{N}_-)=W_{2NF}\pi_{2-}E^2(G_+,\mathfrak{N}_-)\,$. Therefore,
$\,W_{1NF}\pi_{1-}\chi_-^{-1}E_{\alpha}^2(G_+,\mathfrak{N}_-)=
W_{2NF}\pi_{2-}\chi_-^{-1}E_{\alpha}^2(G_+,\mathfrak{N}_-)\,$, where
$\,E_{\alpha}^2(G_+,\mathfrak{N}_-)=\chi_-E^2(G_+,\mathfrak{N}_-)\,$ is the subspace of
character-automorphic functions corresponding to the unitary representation $\,\alpha\,$ of the
fundamental group of the domain $G_+$. In view that
$\,W_{kNF}\pi_{k-}\chi_-^{-1}|E_{\alpha}^2(G_+,\mathfrak{N}_-),\,k=1,2\,$ are isometries, by the
generalized Beurling theorem, there exists an unitary operator
$\,U\in\mathcal{L}(\mathfrak{N}_-)\,$ such that
$\,W_{1NF}\pi_{1-}\chi_-^{-1}=W_{2NF}\pi_{2-}\chi_-^{-1}U\,$. Whence we have
\[
\begin{array}{lcl}
\chi_-\Theta_2^+\chi_+^{-1}&=&\chi_-\pi_{2-}^{\dag}\pi_{2+}\chi_+^{-1}
=(\pi_{2-}\chi_-^{-1})^*\pi_{2+}\chi_+^{-1}=\\[2pt]
&=&(W_{2NF}^{-1}W_{1NF}\pi_{1-}\chi_-^{-1}U^{-1})^*W_{2NF}^{-1}W_{1NF}\pi_{1+}\chi_+^{-1}=\\[2pt]
&=&U(\pi_{1-}\chi_-^{-1})^*\pi_{1+}\chi_+^{-1}=
U\chi_-\pi_{1-}^{\dag}\pi_{1+}\chi_+^{-1}=U\chi_-\Theta_1^+\chi_+^{-1}\,.\quad\;\square
\end{array}
\]
\begin{remark}\label{rem9}
In the proof we have made use of the similarity of the bundle shift
$\,z|E_{\alpha}^2(G_+,\mathfrak{N})\,$ to the trivial shift
$\,z|E^2(G_+,\mathfrak{N})\,$~\cite{AD}. This assertion is equivalent to the fact of triviality of
any analytic vector bundles over multiply connected domains (see~\cite{AD}, which in the scalar
case goes back to~\cite{Sa}). We avoid to deal with multiply valued analytic functions and use the
technique of weights instead of the traditional uniformization technique or analytic vector
bundles. Note that the use of weights $\Xi_{\pm}$ can be effective for the simple connected domains
too (see Remarks~\ref{rem14} and~\ref{rem16}).
\end{remark}

\begin{theorema}
If $\;\Theta_1,\,\Theta_2\in\Cfn\,$ and
$\,\mathcal{F}_{ct}(\Theta_1)=\mathcal{F}_{ct}(\Theta_2)\in\mathcal{N}(G_+\cup
G_-,\mathcal{L}(\mathfrak{N}_{-},\mathfrak{N}_{+}))\,$, then $\;\Theta_1=\Theta_2\,$.
\end{theorema}
\noindent \textit{Proof of Theorem A}.\, Without loss of generality we can assume that the transfer
function $\,\Upsilon=\mathcal{F}_{ct}(\Theta_1)=\mathcal{F}_{ct}(\Theta_2)\,$ is $\Xi$-pure (see
Remark~\ref{rem7}). By Prop.~\ref{delta} and Prop.~\ref{mod_space}, we obtain
$\,\chi_-\Theta_2^+\chi_+^{-1}=U\chi_-\Theta_1^+\chi_+^{-1}\,$. Whence it follows the
one-valuedness of the operator valued function $\,\chi_-^{-1}U\chi_-\,$. Whence we can see that the
operator $U$ commutes with the $W^*$-algebra $\,\mathcal{W}^*(\alpha)\,$ generated by the unitary
representation $\alpha$. By the theorem of Bungart~\cite{Bu}, the analytic vector bundles
corresponding to the operator valued character-automorphic function $\,\chi_-\,$ is trivial with
respect to the group of invertible operators $G(\alpha)$ of the algebra
$\,\mathcal{W}^*(\alpha)\,$. That means there exists character-automorphic function $\chi(z)$ with
values in $G(\alpha)$ such that this function trivialize the analytic vector bundle. Therefore the
operator valued function $\,\eta_-=\chi_-^{-1}\chi\,$ is one-valued and $\,\eta_-,\eta_-^{-1}\in
H^{\infty}(G_+,\mathfrak{N}_-)\,$.

Let $\,\eta=(z, I, \eta_-)\,$, $\,\Theta_{1\eta}=\Phi_{\eta}^{\Cfn}(\Theta_{1})\,$, and
$\,\Theta_{2\eta}=\Phi_{\eta}^{\Cfn}(\Theta_{2})\,$. Then $\,
\mathcal{F}_{tc}(\Theta_{1\eta})=\mathcal{F}_{tc}(\Phi_{\eta}^{\Cfn}(\Theta_{1}))=
\Phi_{\eta}^{\Tfn}(\mathcal{F}_{tc}(\Theta_{1}))=\Phi_{\eta}^{\Tfn}(\mathcal{F}_{tc}(\Theta_{2}))=
\mathcal{F}_{tc}(\Phi_{\eta}^{\Cfn}(\Theta_{2}))=\mathcal{F}_{tc}(\Theta_{2\eta}) \,$. Therefore
the transfer function
$\,\Upsilon_{\eta}=\mathcal{F}_{tc}(\Theta_{1\eta})=\mathcal{F}_{tc}(\Theta_{2\eta})\,$ corresponds
to both characteristic functions $\Theta_{1\eta}$ and $\Theta_{2\eta}$. We have
\[
\begin{array}{lcl}
U\Theta_{1\eta}^+&=&U\eta_-^{-1}\Theta_{1}^+=U\chi^{-1}\chi_-\Theta_{1}^+\chi_+^{-1}\chi_+=\\[3pt]
&=&\chi^{-1}(U\chi_-\Theta_{1}^+\chi_+^{-1})\chi_+=\chi^{-1}\chi_-\Theta_{2}^+\chi_+^{-1}\chi_+=
\eta_-^{-1}\Theta_{2}^+=\Theta_{1\eta}^+\,.
\end{array}
\]
Since $\,\Xi_{1\eta-}=\Xi_{2\eta-}=\chi^*\chi\,$ and $\,U\chi^*\chi=\chi^*\chi U\,$, we also have
\[
\begin{array}{lcl}
\Theta_{1\eta}^-U^*&=&\Xi_+^{-1}\Theta_{1\eta}^{+*}\chi^*\chi U^*=
\Xi_+^{-1}(U\Theta_{1\eta}^{+})^*\chi^*\chi=\\[3pt]
&=&\Xi_+^{-1}\Theta_{2\eta}^{+*}\chi^*\chi =\Theta_{2\eta}^-\,.
\end{array}
\]
Then, taking into account (CtoT), we get
$\,\Upsilon_{2\eta}(z)=\Upsilon_{1\eta}(z)U=\Upsilon_{2\eta}(z)U\,$. Whence,
$\,\Upsilon_{2\eta}(z)(U-I)=0\,$. Since, the transfer function $\,\Upsilon_{2\eta}\,$ is\,
$\Xi$-pure, we get $\,U=I\,$ (see Remark~\ref{rem7}). Thus, $\,\Theta_1^+=\Theta_2^+\,$. This
completes the proof of Theorem A. \proofend
\begin{remark}\label{rem10}
If $\,\dim\mathfrak{N}_-<\infty\,$, we can make use of the more weak assertion from~\cite{Gr}
instead of
 Bungart's theorem~\cite{Bu}.
\end{remark}
Analyzing the above proof we can reveal the restoration procedure of the characteristic function
$\,\Theta\,$ for a given transfer function $\,\Upsilon\,$. In the same way as in the proofs of
Prop.~\ref{rk} and Prop.~\ref{mod_space}, we define $\,\mathcal{H}_{NF+}=\mathcal{K}_{\Theta
NF}\dot+E^2(G_+,\mathfrak{N}_+)\,$, where $\,\mathcal{K}_{\Theta
NF}=\vee_{z\in\rho(\widehat{T})}\hat{r}_{nz}\,$. Clearly,
$\,\mathcal{H}_{NF+}\subset\mathcal{H}_{NF}\,$.

Let $\,\mathcal{V}=\vartheta(\mathcal{U})\,$, where $\,\mathcal{U}\,$ is the operator of
multiplication by the independent variable in $\,\mathcal{H}_{NF}\,$ and $\,\vartheta\,$ is a
scalar inner function on $\,G_+\,$ with precisely $n$ zeroes (we assume that the domain $\,G_+\,$
has $\,n\,$ boundary components)~\cite{Ba}. Since the subspace $\,\mathcal{H}_{NF+}\,$ is invariant
under the operator $\,\mathcal{V}\,$, the operator $\,\mathcal{V}|\mathcal{H}_{NF+}\,$ is an
isometry. Then there exists the Wold decomposition of the subspace $\,\mathcal{H}_{NF+}\,$ with
regard to the operator $\,\mathcal{V}|\mathcal{H}_{NF+}\,$ (see~\cite{NF}). We put
$\,\mathfrak{M}=\mathcal{H}_{NF+}\ominus\mathcal{V}\mathcal{H}_{NF+}\,$. Taking into account the
vector representation of Hardy spaces for multiply connected domains (see, e.g.,~\cite{F}), we see
$\,\dim\mathfrak{M}=n\cdot\dim\mathfrak{N}_-\,$. Let
$\,V\in\mathcal{L}(\mathfrak{N}_-^n,\mathfrak{M})\,$ be a unitary operator. Define
\[
\begin{array}{ll}
Wf:=\sum\limits_{k=-\infty}^{+\infty} \mathcal{U}_0^k
V^{-1}P_{\mathfrak{M}}\mathcal{V}^{-k}f\,,\quad  & W : \mathcal{H}_{NF} \to
L^2(\mathbb{T},\mathfrak{N}_-^n)\,,
\end{array}
\]
where $\,\mathcal{U}_0\,$ is the operator of multiplication by the independent variable in
$\,L^2(C,\mathfrak{N}_-)\,$.

Similarly, let $\,\mathcal{V}_0=\vartheta(z)\,$ in the space $\,L^2(C,\mathfrak{N}_-)\,$ and
$\,\mathfrak{M}_0=E_{\alpha}^2(G_+,\mathfrak{N}_-)\ominus\mathcal{V}_0E_{\alpha}^2(G_+,\mathfrak{N}_-)\,$,
where the unitary representation $\alpha$ is corresponded to the outer operator valued function
$\,\chi_-\,$ defined by $\,\chi_-^*\chi_-=\Xi_-\,$. Let
$\,V_0\in\mathcal{L}(\mathfrak{N}_-^n,\mathfrak{M}_0)\,$ be a unitary operator. It can be shown
that
\[
L^2(C,\mathfrak{N}_-)=\bigoplus_{k=-\infty}^{+\infty}\,\mathcal{V}_0^k \mathfrak{M}_0\,,\quad
E_{\alpha}^2(G_+,\mathfrak{N}_-)=\bigoplus_{k=0}^{+\infty}\,\mathcal{V}_0^k \mathfrak{M}_0\,.
\]
We define also
\[
\begin{array}{ll}
W_0f:=\sum\limits_{k=-\infty}^{+\infty} \mathcal{U}_0^k
V_0^{-1}P_{\mathfrak{M}_0}\mathcal{V}_0^{-k}f\,,\quad  & W_0 : L^2(C,\mathfrak{N}_-) \to
L^2(\mathbb{T},\mathfrak{N}_-^n)\,.
\end{array}
\]
It is easily shown that $\,W^*\,$ is an isometry, $\,W_0\,$ is an unitary operator, and
$\,W\mathcal{V}=\mathcal{U}_0W\,,W_0\mathcal{V}_0=\mathcal{U}_0W_0\,$. Besides, we have
$\,W_0E_{\alpha}^2(G_+,\mathfrak{N}_-)=H^2(\mathbb{D},\mathfrak{N}_-^n)\,$ and
$\,W\mathcal{H}_{NF+}^p=H^2(\mathbb{D},\mathfrak{N}_-^n)\,$, where
$\,\mathcal{H}_{NF+}^p:=\oplus_{k=0}^{+\infty}\,\mathcal{V}^k \mathfrak{M}\,$ and the subspace
$\,\mathcal{H}_{NF+}^u:=\mathcal{H}_{NF+}\ominus\mathcal{H}_{NF+}^p\,$ reduces the operator
$\,\mathcal{U}\,$.

\vskip 2pt  Let an unitary operator $\,Y : \mathfrak{N}_-^n\to\mathfrak{N}_-^n\,$ be a solution of
the linear equation $\,\mathcal{U}W^*m_YW_0|\mathfrak{M}_0=W^*m_YW_0z|\mathfrak{M}_0\,$, where
$\,m_Y\,$ is the operator of multiplication by the unitary constant $\,Y\,$ in
$\,L^2(\mathbb{T},\mathfrak{N}_-^n)\,$. Such a solution necessarily exists if $\,\Upsilon(z)\,$ is
a transfer function. Note that the equation is equivalent to the Riccati equation
$\,Y^*W\mathcal{U}W^*Y|\mathfrak{N}_-^n=W_0zW_0^*|\mathfrak{N}_-^n\,$. If $\,Y\,$ is a solution of
the equation, the identity $\,\mathcal{U}W^*m_YW_0=W^*m_YW_0z\,$ can easily be extended to the
whole space $\,L^2(\mathbb{T},\mathfrak{N}_-^n)\,$.

We set $\,\pi_{0+}:=\left(\begin{array}{l}I\\0\end{array}\right) : L^2(C,\mathfrak{N}_+)
\to\mathcal{H}_{NF}\,$ and $\,\pi_{0-}:=W^*m_YW_0 : L^2(C,\mathfrak{N}_-) \to\mathcal{H}_{NF}\,$.
It can easily be checked that $\,\pi_{0-}^*\pi_{0-}=I\,$, $\,\mathcal{U}\pi_{0-}=\pi_{0-}z\,$, and
$\,\pi_{0-}E_{\alpha}^2(G_+,\mathfrak{N}_-)=\mathcal{H}_{NF+}^p\,$. A mapping $\,\pi_{0-}'\,$
satisfying these conditions is determined up to a constant unitary operator $\,U\,$ such that
$\,\pi_{0-}'=\pi_{0-}U^*\,$.

Let $\,\Theta_U^{+}:=\chi_-^{-1}U\Theta_0^{+}\chi_+\,$ and
$\,\Theta_U^{-}:=(\chi_+^*\chi_+)^{-1}\Theta_U^{+*}(\chi_-^*\chi_-)\,$, where
$\,\Theta_0^{+}=\pi_{0-}^*\pi_{0+}\,$. According to Theorem A the unitary operator $\,U\,$ is
uniquely determined by the conditions
$\,\Upsilon_+(\zeta)=\Theta^{-}_{U+}(\zeta)-\Theta_{U}^+(\zeta)^{-1},\;
\,\Upsilon_-(\zeta)=-\Theta^-_{U-}(\zeta)\,,\;\zeta\in C\,$. Note that these conditions are linear.
Thus we can complete recover the characteristic function $\Theta$ for a given transfer function
$\Upsilon$.
\begin{remark}\label{rem11}
For simple-connected domains, we can suggest a more direct procedure~\cite{T3}, which essentially
exploits the relation $\rm{(\Phi\&F)}$ and reduces the problem to the case of the unit disk and
unitary colligations.
\end{remark}
\begin{remark}\label{rem12}
If the absolutely continuous spectrum of the operator $\,T\,$ is rich , e.g.,
$\,\clos\Delta^+(\zeta)\,\mathfrak{N}_-=\mathfrak{N}_-\,$, then, by Lemma~\ref{l_delta}, we can
restore the characteristic function $\Theta^+$ from the relation
$\,\Theta^{-}_{+}(\zeta)-\Theta^+(\zeta)^{-1}=\Upsilon_+(\zeta)-\Upsilon_-(\zeta)\,,\;\zeta\in
C\,$.
\end{remark}
\begin{remark}\label{rem13}
We have established the Theorem A under the assumption that the weights $\,\Xi_1,\,\Xi_2\,$ are
equal. If this assumption is dropped, we lose the uniqueness and the characteristic functions
$\Theta_1^+,\Theta_2^+$ satisfy some relations (see~\cite{Ya}). Note that in his paper D.Yakubovich
considered only the case of $\Xi$-inner operator valued functions, i.e., he assumed (in our terms)
that $\,\Theta^-=(\Theta^+)^{-1}\,$.
\end{remark}

We pass to the problem to construct the functional model for a given simple (see Remark~\ref{rem4})
conservative curved system. Let $\,\Sigma=(T,M,N,\Theta_u,\Xi)\,$. So, we can consider its transfer
function $\Upsilon$. If $\,\Upsilon(z)\in\mathcal{N}(G_+\cup
G_-,\mathcal{L}(\mathfrak{N}_{-},\mathfrak{N}_{+}))\,$, then by the above we are able to restore
the corresponding characteristic function $\Theta\,$. Thus we need to verify that certain
construction (=the functional model) agrees with the system $\,\Sigma\,$. Such a kind of problem
arises in applications~\cite{VYa} and was studied in~\cite{Ya} for the case of $\Xi$-inner
functions.

The operator $\,W\,$, which realizes a similarity of a system $\,\Sigma\,$ to the model system
$\,\widehat{\Sigma}\,$, can easily be calculated: $\,W=W_{NF}^{-1}V\,$, where
$\,V\in\mathcal{L}(H,\mathcal{K}_{NF})\,$ and $\,\forall f\in H$
\[
\begin{array}{l}
(Vf)_{\pi}:=-[M(T-\zeta)^{-1}f]_-\,,\\[6pt]
(Vf)_{\tau}:=(\Delta^+)^{-1}\Theta^+([M(T-\zeta)^{-1}f]_+-[M(T-\zeta)^{-1}f]_-)\,.
\end{array}
\eqno{\rm{(Sim)}}
\]
Therefore we need to verify accordance of the operator $\,V\,$ and the functional model.
\begin{prop}\label{StoM}
Suppose a system $\,\Sigma=(T,M,N,\Theta_u,\Xi)\,$ is simple and there exists $\,\Theta\in\Cfn\,$
such that $\,\Upsilon=\mathcal{F}_{ct}(\Theta)\,$ and $\,\Theta^+(z)^{-1}$ possesses boundary
values a.e. on $C$, where $\,\Upsilon(z)=M(T-z)^{-1}N\,$. Suppose that $\,\forall f\in H\,$ there
exist ${\,[M(T-\zeta)^{-1}f]_{\pm}\,}$ a.e. on the curve $C$. Let $\,V : H\to\mathcal{K}_{NF}\,$ be
defined by formulas {\rm(Sim)}. If $\,||V||<\infty\,$ and $\,||V^{-1}||<\infty\,$. Then
$\,\Sigma\in\Sys\,$.
\end{prop}
\myproof It follows from the identity $\,V{r}_{nz}=\hat{r}_{nz}\,$. \proofend More symmetric
variant of this assertion can be formulated if we take into account the dual model.
\begin{prop}
Under the hypotheses of Proposition~\ref{StoM} let the operator $\,V_* : H\to\mathcal{K}_{*NF}\,$
be defined by formulas those dual to {\rm(Sim)}. If $\,||V||<\infty\,$ and $\,||V_{*}||<\infty\,$.
Then $\,\Sigma\in\Sys\,$.
\end{prop}
\myproof Additionally to the identity $\,V{r}_{nz}=\hat{r}_{nz}\,$ we need to make use of the
identity $\,(V{r}_{nz},V_*{r}_{*mw})=({r}_{nz},{r}_{*mw})\,$. \proofend
\begin{remark}\label{rem14}
For simple-connected domains and scalar weights, we can present another description for
conservative curved systems~\cite{T1,T2,T3} :
\[
(T,M,N)\;=\;(W_0\varphi(T_0)W_0^{-1},\,M_0\psi_+(T_0)W_0^{-1},\,W_0\psi_-(T_0)N_0)\,,
\]
where $\,{\mathfrak A}_0=\left(%
\begin{array}{cc}
  T_0 & N_0 \\
  M_0 & L_0 \\
\end{array}%
\right)\,$ is a simple unitary colligation, $\,W_0\in\mathcal{L}(H_0,H)\,$,
$\,W_0^{-1}\in\mathcal{L}(H,H_0)\,$, $\,\varphi : \mathbb{D}\to G_{+}\,$ is a conformal map,
$\,\psi_+=\sqrt{\varphi '}/(\eta_+\circ\varphi)\,,\;\, \psi_-=\sqrt{\varphi
'}\,(\eta_-\circ\varphi)\,$, and $\,\eta_{\pm},\eta_{\pm}^{-1}\in
H^{\infty}(G_{+},\mathcal{L}(\mathfrak{N}_{\pm}))\,$.
\end{remark}
\begin{remark}\label{rem15}
Note also that  $\,T\,$ is the main operator of a system $\,\Sigma\in\Sys\,$ iff $\,T\,$ is similar
to a c.n.u. $G_+$-contraction. An operator $\,T_0\in\mathcal{L}(\mathcal{K}_0)\,$ is called a
$G_+$-contraction if there exists a Hilbert space $\,\mathcal{H}_0\,$ and a normal operator
$\,\mathcal{U}_0\in\mathcal{L}(\mathcal{H}_0)\,$ such that $\,\mathcal{K}_0\subset\mathcal{H}_0\,$,
$\,(T_0-z)^{-1}=P_0(\mathcal{U}_0-z)^{-1}|\mathcal{K}_0\,,\;z\in G_-\,$,
$\,\sigma(\mathcal{U}_0)\subset C\,$, where $\,P_0\,$ is the orthogonal projection onto the
subspace $\,\mathcal{K}_0\,$. A $G_+$-contraction $\,T_0\,$ is called a c.n.u. $G_+$-contraction if
the operator $\,T_0\,$ have no nontrivial subspaces $\,\mathcal{K}_0'\subset\mathcal{K}_0\,$ such
that the operator $\,T_0|\mathcal{K}_0'\,$ is normal and $\,\sigma(T_0|\mathcal{K}_0')\subset C\,$.
\end{remark}

\section{Applications to perturbation theory}

In this Section we employ the functional model to establish duality of spectral
components~\cite{T1}. Note that the duality problem was parental for our variant of generalized
Sz.-Nagy-Foia\c{s}'s model.

Recall the definitions of the spectral components~\cite{Na2,T1,MV,Ma}:
\[
\begin{array}{lcl}
{\widetilde M}(S)&=&\{\;f\in H :\;\forall g\in H\quad
((S-\zeta)^{-1}f,g)_+=((S-\zeta)^{-1}f,g)_-\,,\; \zeta\in C\,\}\,;\\
{\widetilde N}_{\pm}(S)&=&\{\;f\in H :\;\forall g\in H\quad
((S-z)^{-1}f,g)_{\pm}\in E^2(G_{\pm})\}\,;\\
{\widetilde D}_{\pm}(S)&=&\{\;f\in H :\;\forall g\in H\quad ((S-z)^{-1}f,g)_{\pm}\in
D(G_{\pm})\}\,,
\end{array}
\]
where $D(G_{\pm})=\{\, f : f=\varphi/\delta,\;\, \varphi,\delta\in H^{\infty}(G_{\pm}),
\;\,\delta\,\;\mbox{is an outer function}\}$ is the Nevanlinna-Smirnov space~\cite{Du}. Various
combinations of the above linear subspaces are considered: $\, {{\widetilde {N}}={\widetilde
N}_+\cap {\widetilde N}_-},\; {{\widetilde {NM}}_{\pm}={\widetilde N}_{\pm}\cap {\widetilde M}},\;
{{\widetilde {DM}}_{\pm}={\widetilde D}_{\pm}\cap {\widetilde M}}\,$. The closure ${X(S)={\rm
clos\,}{\widetilde X}(S)}\,$, where ${\widetilde X}\in\{{\widetilde M},{\widetilde D}_{\pm},
{\widetilde {DM}}_{\pm},{\widetilde N},{\widetilde N}_{\pm}, {\widetilde {NM}}_{\pm}\}$, is called
a (weak) spectral component of the operator $S$. Note that the absolutely continuous $\,N(S)\,$ and
singular $\,M(S)\,$ subspaces are generalizations of the corresponding notions for self-adjoint
operators~\cite{Yaf} and for contractions of $C_{11},\,C_{00}$ classes~\cite{NF}. Concerning of a
spectral meaning of these components we refer the reader to~\cite{Na2,T1,Ma}. The main theorem
from~\cite{T1} is the following.
\begin{theoremb}\label{dual}
Suppose $C$ is a simple closed $C^{4+\varepsilon}$-smooth curve and  $\,U,S\in \mathcal{L}(H)\,$
are operators such that $\,U^*U=UU^*,\, S-U\in \mathfrak{S}_1,\,\sigma(U)\subset C\,,
\sigma_c(S)\subset C$.\, Then
\[
1)\;N(S)^\bot=M(S^*)\,;,\quad 2)\;N_{\pm}(S)^\bot=DM_{\mp}(S^*)\,;\quad
3)\;NM_{\pm}(S)^\bot=D_{\mp}(S^*)\,.
\]
\end{theoremb}
\noindent In the present paper we extend this result to the case of multiply connected domains.
Note that particular cases of this Theorem play an important role in non-selfadjoint scattering
theory and in the theory of extreme factorizations of J -contraction-valued functions (J
-outer-inner and A-singular-regular factorizations) (see the corresponding references
in~\cite{T1}).

As in~\cite{T1} we are going to solve the duality problem by means of the functional model. Since
the spectral components appear in the duality relations in a symmetrical way, the both Smirnov
classes $E^2(G_+)$ and $E^2(G_-)$ ought to be tantamount ingredients of the functional model. By
the same reason, the standard trick to involve into a construction the Riemann surface (= double of
$G_+$, see, e.g.,~\cite{PF,F}) does not work in our case because we lose all geometric information
concerning the domain~$G_-$.

We shall use generalized model of S.Naboko~\cite{Na,Na2,MV}. The key observation was done
in~\cite{Na} (see also~\cite{Na2,MV}), where it was revealed the possibility to study effectively
by means of Sz.-Nagy-Foia\c{s}'s model perturbations of certain form. It is obvious that any trace
class perturbation of a self-adjoint operator can be represented in the form
$\,L=(L_R+i|L_I|)+|L_I|^{1/2}\varkappa |L_I|^{1/2}\,$, where $\,L_R=\frac{1}{2}
(L+L^*),\,L_I=\frac{1}{2i} (L-L^*)\,$ and $|L_I|^{1/2}\in\mathfrak{S}_2$. It is an
exercise~\cite{MV} to verify that any trace class perturbation of an unitary operator can be
represented in the form $\,S=T+(I-TT^*)^{1/2}\varkappa(I-T^*T)^{1/2}\,$, where $T$ is a contraction
and $\,(I-TT^*)^{1/2},\,(I-T^*T)^{1/2}\in\mathfrak{S}_2\,$. Therefore any trace class perturbation
of self-adjoint or unitary operator can be studied using Naboko's approach. In our setting these
perturbations look like $\,{S}={T}+{N}\varkappa{M}\,$, where $\,({T},{M},{N})\,$ is a conservative
curved system and $\varkappa\,$ is a bounded linear operator. It is a non-trivial problem to prove
that any trace class perturbation of a normal with the spectrum lying on a curve $C$ admits such
representation. It was the main obstacle, which we surmounted in~\cite{T1}. For
$C^{4+\varepsilon}$-smooth curves we have established the following theorem.
\begin{theorem_}\label{repr}
Suppose $C$ is a simple closed $C^{4+\varepsilon}$-smooth curve and  $\,U,S\in \mathcal{L}(H)\,$
are operators such that $\,U^*U=UU^*,\, S-U\in \mathfrak{S}_1,\,\sigma(U)\subset C\,,
\sigma_c(S)\subset C$.\, Then the operator $S$ admits the representation
$S=\varphi(T_0)+N_0\varkappa M_0$, where $\varphi$ is a conformal mapping of the unit disk ${\Bbb
D}$ onto the domain $G_+$,
$\mathfrak{A}_0=\left(\begin{array}{cc}T_0&N_0\\M_0&L_0\end{array}\right)$ is a simple unitary
colligation, $T_0^{-1}\in\mathcal{L}(H),\; M_0,N_0\in \mathfrak{S}_2,\;
\varkappa\in\mathcal{L}(\mathfrak{N})$ (without loss of generality it can be assumed that
$\,\mathfrak{N}_+=\mathfrak{N}_-=\mathfrak{N}\,$).
\end{theorem_}
\begin{remark}\label{rem16}
In this connection note that this representation lead us to the model with the scalar weights
$\,\Xi_+=|\varphi'| I,\;\Xi_-=\frac{1}{|\varphi'|} I\,$ (cf. Remark~\ref{rem9}).
\end{remark}
\noindent In~\cite{T1} the theorem was proved under assumption that the domain $G_+$ is
simple-connected. Actually, it can be extended to the case of multiply connected domains.
\begin{prop}\label{mc_repr}
Suppose $C$ is a simple closed $C^{4+\varepsilon}$-smooth curve and  $\,U,S\in \mathcal{L}(H)\,$
are operators such that $\,U^*U=UU^*,\, S-U\in \mathfrak{S}_1,\,\sigma(U)\subset C\,,
\sigma_c(S)\subset C$.\, Then there exists $\,\Pi\in\Mod\,$ and an invertible operator $W$ such
that $\,S=W^{-1}\widehat{S}W=W^{-1}(\widehat{T}+\widehat{N}\varkappa\widehat{M})W\,$,
$\,(\widehat{T},\widehat{M},\widehat{N})=\mathcal{F}_{sm}(\Pi)\,$,
$\,\widehat{M},\widehat{N}\in\mathfrak{S}_2,\,\varkappa\in\mathcal{L}(\mathfrak{N})\,$.
\end{prop}
\noindent \textit{Sketch of the proof.}\,\, Let $\,C=(\cup_{k=1}^n C_k)\cup C_{0}\,$,
$\,G_-=(\cup_{k=1}^n G_{k-})\cup G_{0-}\,$, where $\,G_{k-}:=\Int C_k\,$ and $\,\infty\in
G_{0-}:=\Ext C_0\,$. Using spectral projections from Riesz-Danford calculus, we can construct the
operator $\,\oplus_{k=0}^n S_k\,$ such that $\,S_k-U_k\in \mathfrak{S}_1,\,U_k^*U_k=U_kU_k^*,\,
\sigma(U_k)\subset C_k\,, \sigma_c(S_k)\subset C_k$  and the operator $\,\oplus_{k=0}^n S_k\,$ is
similar to the operator $S$. Note that $U_k$ are not necessary the spectral parts of the operator
$U$ corresponding to the components $C_k$ (naive candidates to realize similarity can have nonzero
index; recall that $\,\ind V:=\dim\Ran V - \dim\Ker V\,$). In fact, $U_k$ may be obtained using
those spectral parts by adding or removing eigenvectors to the corresponding subspace. Then we can
reduce our problem to the simple-connected case: $\,S_k=T_k+N_k\varkappa M_k\,$. We need only to
note that the models of the operators $\,T_k\,$ corresponding to the curves $C_k$ and $C$ are
similar. \proofend

\vskip -10pt Thus we have reduced the duality problem to operators of the form
$\,\widehat{S}=\widehat{T}+\widehat{N}\varkappa\widehat{M}\,$. Note that, since Naboko's model is a
superstructure on Sz.-Nagy-Foia\c{s}'s model, we made an extra effort in Section~1 (as well as
in~\cite{T1,T2,T3}) to develop our variant of generalized Sz.-Nagy-Foia\c{s}'s model as simple as
possible with the aim to focus entirely upon the superstructure by automating calculations on the
level of the underlying model.

The structure of the perturbation $\widehat{N}\varkappa\widehat{M}$ permits us to compute the
resolvent of the operator $\widehat{S}$ :
\[
(\widehat{S}-z)^{-1}f=(\mathcal{U}-z)^{-1}(f+(\pi_+\varkappa_+^r+\pi_-\varkappa_-^r)\,n(z)),\;\;\;
f\in\mathcal{K}_{\Theta},
\]
where $\quad\,n(z)=\Theta_{\cdot\varkappa}^{\pm}(z)^{-1}(\pi_{\mp}^{\dag}f)(z)\,,\quad
\varkappa_+^r=-I-\Theta_+^-\varkappa\,,\quad \varkappa_-^r=\varkappa\,,\,$\\[2pt]
$\,\Theta_{\cdot\varkappa}^{+}=-\varkappa+\Theta^++\Theta^+\Theta_+^-\varkappa\,,\quad
\Theta_{\cdot\varkappa}^{-}=I-\Theta_-^-\varkappa\,$\,. Using this formula we obtain the
descriptions of the spectral components in terms of the functional model:
\[
\begin{array}{lcl}
\widetilde{M}(\widehat{S})&=&\{f\in\mathcal{K}_{\Theta} \,:\,
\Theta_{\cdot\varkappa}^{-}(\zeta)^{-1}(\pi_{+}^{\dag}f)(\zeta)=
\Theta_{\cdot\varkappa}^{+}(\zeta)^{-1}(\pi_{-}^{\dag}f)(\zeta)\, \}\,,\\[5pt]
\widetilde{N}_{\pm}(\widehat{S})&=&\{f\in\mathcal{K}_{\Theta} \,:
\pi_{\mp}^{\dag}f\in\Theta_{\cdot\varkappa}^{\pm}(z)\,E^2(G_{\pm},\mathfrak{N})\, \}\,,\\[5pt]
\widetilde{D}_{\pm}(\widehat{S})&=&\{f\in\mathcal{K}_{\Theta} \,:\,
\pi_{\mp}^{\dag}f\in\Theta_{\cdot\varkappa i}^{\pm}(z)\,E^2(G_{\pm},\mathfrak{N})\,\}\,,
\end{array}
\]
where $\,\Theta_{\cdot\varkappa}^{\pm}(z)=\Theta_{\cdot\varkappa
i}^{\pm}(z)\,\Theta_{\cdot\varkappa e}^{\pm}(z)\,$ is the inner-outer factorization~\cite{NF}.
Since there exists operators $W$ and $U$ such that $\,W^{-1}\widehat{S}W-U\in\mathfrak{S}_1\,$,
$\,U^*U=UU^*,\,\sigma(U)\subset C\,$ and $\,\widehat{M},\widehat{N}\in\mathfrak{S}_2\,$, we have
$\,(\Theta_{\cdot\varkappa}^{\pm})^{-1}\in\mathcal{N}(G_{\pm},\mathcal{L}(\mathfrak{N}))$. Besides,
if the curve $C$ is $C^{2+\varepsilon}$ smooth then $\;\Theta_{\pm}^-\in
H^{\infty}(G_{\pm},\mathcal{L}(\mathfrak{N}))\,$. We "lift" spectral components up to the level of
the "dilation" space $\mathcal{H}$\,:
\[
\begin{array}{lcl}
M(\Pi,\varkappa) &=& \{\,f\in \mathcal{H} :\,
\Theta_{\cdot\varkappa}^{-}(\zeta)^{-1}(\pi_{+}^{\dag}f)(\zeta)=
\Theta_{\cdot\varkappa}^{+}(\zeta)^{-1}(\pi_{-}^{\dag}f)(\zeta)\,\}\,,\\[5pt]
N_{\pm}(\Pi,\varkappa) &=& \{\,f\in \mathcal{H}
:\,P_{\pm}(\varkappa_+^l\pi_-^{\dag}+\varkappa_-+^l\pi_+^{\dag})f=0 \,\}\,,\\[5pt]
D_{\pm}(\Pi,\varkappa) &=& \{\,f\in \mathcal{H} :\,\pi_{\mp}^{\dag}f\in\Theta_{\cdot\varkappa
i}^{\pm}\,E^2(G_{\pm},\mathfrak{N}) \,\}\,.
\end{array}
\]
For them we have
\[
\widetilde{M}(\widehat{S})=M(\Pi,\varkappa)\cap\mathcal{K}_{\Theta}\,,\;\;
\widetilde{N}_{\pm}(\widehat{S})=P_{\Theta}N_{\pm}(\Pi,\varkappa)\,,\;\;
\widetilde{D}_{\pm}(\widehat{S})=D_{\pm}(\Pi,\varkappa)\cap\mathcal{K}_{\Theta}\,.
\]
The possibility itself of such lifting depends heavily on the structure of the perturbation in
question. The above mentioned properties
$\,(\Theta_{\cdot\varkappa}^{\pm})^{-1}\in\mathcal{N}(G_{\pm},\mathcal{L}(\mathfrak{N}))$ and
$\;\Theta_{\pm}^-\in H^{\infty}(G_{\pm},\mathcal{L}(\mathfrak{N}))\,$ are also crucial for this
lifting.

Using the dual model one can prove that $\,N(\Pi,\varkappa)^{\bot}=M(\Pi_*,\varkappa_*)\,$,
$\,N_{\pm}(\Pi,\varkappa)^{\bot}=DM_{\mp}(\Pi_*,\varkappa_*)\,$, and
$\,NM_{\pm}(\Pi,\varkappa)^{\bot}=D_{\mp}(\Pi_*,\varkappa_*)\,$. Then we can rewrite the proof
from~\cite{T1} almost literally. We need only to use some notions from~\cite{T1} in extended
meaning, e.g., in Prop.4.9(4)~\cite{T1} we need to replace the subspaces
$\,E_{\Theta_{\varkappa\cdot *i}^{\pm}}\,$ by $\,E_{\Omega^{\pm}\Theta_{\varkappa\cdot
*i}^{\pm}}\,$, where $\,\Omega^{\pm}\,$ is multipliers that countervails a possible
multiple-valuedness of $\,\Theta_{\varkappa\cdot *i}^{\pm}\,$. If the boundary of the domain $G_+$
is $C^{2+\varepsilon}$ smooth, the operator valued functions $\,\Omega^{\pm}\,$ can be chosen
$C^{2+\varepsilon}$ smooth too. The latter has as consequence~\cite{T1,T3} that the weak smooth
vectors are dense in the subspaces $\,{N}_{\pm}(\widehat{S})\,$.

\end{document}